\documentclass[12pt]{amsart} 
\usepackage{amscd, amssymb}
\theoremstyle{plain}
\newtheorem{proposition}{Proposition}
\newtheorem{lemma}{Lemma}
\newtheorem{theorem}{Theorem}
\newtheorem{corollary}{Corollary}

\theoremstyle{definition}

\theoremstyle{remark}

\newcommand{\pp}{\mathbb{P}}
\newcommand{\zz}{\mathbb{Z}}
\newcommand{\qq}{\mathbb{Q}}
\newcommand{\cc}{\mathbb{C}}

\newcommand{\mbar}{{\overline{M}}}
\newcommand{\so}{{\mathcal{O}}}
\newcommand{\Pic}{{\mathcal{P}}ic}
\newcommand{\cl}{\mathcal{L}}

\newcommand{\ev}{{\operatorname{ev}}}

\newcommand{\vir}{\operatorname{vir}}
\newcommand{\la}{\langle}
\newcommand{\ra}{\rangle}

\begin{document}
\title[A Reconstruction theorem]{A reconstruction
theorem in quantum cohomology and quantum $K$-theory}
\author{Y.-P.~Lee} 
\address{UCLA}
\email{yplee@math.ucla.edu}
\author{R.~Pandharipande}
\address{CalTech}
\email{rahulp@its.caltech.edu}
\thanks{Y.P. was partially supported by DMS-0072547.  
R.P. was partially supported by DMS-0071473 and
fellowships from the Sloan and Packard foundations.}

\date{}

\begin{abstract}
A reconstruction theorem for genus 0 gravitational
quantum cohomology and quantum $K$-theory is proved.
A new linear equivalence in the Picard group of the
moduli space of genus 0 stable maps relating the pull-backs of line
bundles from the target via different markings is 
used for the reconstruction result. Examples of calculations in
quantum cohomology and quantum $K$-theory are given.

\end{abstract}

\maketitle

\setcounter{section}{-1}
\section{Introduction}
\subsection{Divisor relations}
Let $X$ be a nonsingular, projective, complex algebraic variety. Let 
$L$ be a line bundle on $X$.
The goal of the present article is to study the relationship between
the different evaluation pull-backs
$\ev^*_i(L)$ and $\ev_j^*(L)$ on the space of stable
maps $\mbar_{0,n}(X,\beta)$. We will also examine the relationship between
the cotangent line classes $\psi_i$ and $\psi_j$
at distinct markings.

Our method is to study the relationship first in the case of
projective space.
The moduli of stable maps $\mbar_{0,n}(\pp^r,\beta)$ is
a nonsingular Deligne--Mumford stack.
The Picard group
$\Pic(\mbar_{0,n}(\pp^r,\beta))$  with
$\qq$-coefficients
has been analyzed in \cite{rp}. In case $n\geq 1$, the Picard group
is generated by the evaluation pull-backs
$\ev^*_i(L)$ together with the
boundary divisors.
The cotangent line classes $\psi_i$
%$\cl_i$ (or the corresponding divisor $\psi_i=c_1(\cl_i)$)
also determine elements of the Picard group.

\begin{theorem} \label{t:1}
The following relations hold in $\Pic(\mbar_{0,n}(\pp^r,\beta))$ for all
$L \in \Pic(\pp^r)$ and markings $i\neq j$:
\begin{equation} \label{e:1}
\ev_i^*(L)=
 \ev_j^*(L) + \langle \beta, L \rangle \psi_j 
 - \sum_{\beta_1+\beta_2=\beta} \la\beta_1, L\ra D_{i,\beta_1| j,\beta_2} ,
\end{equation}
\begin{equation} \label{e:2}
 \psi_i +\psi_j = D_{i| j},
\end{equation}
where 
$\langle \beta, L \rangle$ denotes the intersection pairing
$$\langle \beta, L \rangle= \int_{\beta} c_1(L).$$
\end{theorem}
 
The boundary notation used in the Theorem is defined by the
following conventions. Let
 $D_{S_1,\beta_1 | S_2, \beta_2}$
denote the  boundary
divisor in $\mbar_{0,n}(\pp^r,\beta)$ parameterizing maps
with reducible domains and splitting types
$$S_1\cup S_2= \{1,\ldots,n\}, \ \ \beta_1+\beta_2=\beta$$
of the marking set and the degree respectively (see \cite{fp, rp}).
We then define:
$$D_{i,\beta_1|j,\beta_2} = \sum_{i\in S_1,\ j\in S_2} D_{S_1, \beta_1| S_2, \beta_2},$$
$$D_{i|j} = \sum_{i\in S_1,\ j\in S_2,\ \beta_1+\beta_2=\beta} D_{S_1, \beta_1| S_2, \beta_2}.$$

Let $L \in \Pic(X)$ be a very ample line bundle on $X$. Let
$$\iota: X \rightarrow \pp^r$$
be the embedding determined by $L$.
There is a canonically induced embedding:
$$\overline{\iota}: \mbar_{0,n}(X, \beta) \rightarrow \mbar_{0,n}(\pp^r, \iota_*[\beta]).$$
The $\overline{\iota}$ pull-backs of the relations of Theorem 1 together with 
the splitting axiom of Gromov--Witten theory yield relations in
the rational Chow group of $\mbar_{0,n}(X,\beta)$:
\begin{eqnarray*} 
\ev_i^*(L) \cap [\mbar_{0,n}(X,\beta)]^{\vir} & = & 
 (\ev_j^*(L)   + \langle \beta, L \rangle \psi_j)  \cap [\mbar_{0,n}(X,\beta)]^{\vir}  \\ & & 
 - \sum_{\beta_1+\beta_2=\beta} \la\beta_1, L\ra [D_{i,\beta_1| j,\beta_2}]^{\vir} . 
\end{eqnarray*}
\begin{eqnarray*} 
 (\psi_i +\psi_j) \cap [\overline{M}_{0,n}(X,\beta)]^{\vir} & = & [D_{i| j}]^{\vir}.
\end{eqnarray*}
Here, $[D_{i,\beta_1| j,\beta_2}]^{\vir}$ and $[D_{i|j}]^{\vir}$ denote the push-forward
to $\mbar_{0,n}(X,\beta)$ of the virtual classes of their constituent boundary divisors.

We note the classes of very ample divisors span $\Pic(X)$ with
$\qq$-coefficients for projective $X$.
As the above Chow relation is {\em linear} in $L$, we conclude:

\begin{corollary} \label{c:1}
The following relations hold in $A_*(\mbar_{0,n}(X,\beta))$ for all
$L \in Pic(X)$ and markings $i\neq j$:

\begin{eqnarray*} 
\ev_i^*(L) \cap [\mbar_{0,n}(X,\beta)]^{\vir} & = & 
 (\ev_j^*(L)   + \langle \beta, L \rangle \psi_j)  \cap [\mbar_{0,n}(X,\beta)]^{\vir}  \\ & & 
 - \sum_{\beta_1+\beta_2=\beta} \la \beta_1, L\ra [D_{i,\beta_1| j,\beta_2}]^{\vir} , 
\end{eqnarray*}

\begin{eqnarray*} 
 (\psi_i +\psi_j) \cap [\overline{M}_{0,n}(X,\beta)]^{\vir} & = & [D_{i| j}]^{\vir}.
\end{eqnarray*}
\end{corollary}

\subsection{Reconstruction}
We use the following standard notation for the Gromov--Witten 
invariants:
\[
 (\tau_{k_1}(\gamma_1),\cdots, \tau_{k_{n}}(\gamma_n))_{0,n, \beta}
 = \int_{ [\mbar_{0,n}(X,\beta)]^{\vir}} \prod_i \psi_i^{k_i} \ev_i^*(\gamma_i)
\]
where $\gamma_i \in H^*(X)$. The
\emph{quantum $K$-invariants} \cite{yl2} are:
\[
 (\tau_{k_1}(\gamma_1),\cdots, \tau_{k_{n}}(\gamma_n))^K_{0,n, \beta}
 = \chi(\mbar_{0,n}(X,\beta), [\mathcal{O}^{\vir}_{\mbar_{0,n}(X,\beta)}]
  \prod_i [\mathcal{L}_i]^{k_i}  \ev_i^*(\gamma_i) ),
\]
where $\gamma_i \in K^*(X)$, $[\mathcal{O}^{vir}_{\mbar_{0,n}(X,\beta)}] \in 
K_0(\mbar_{0,n}(X,\beta))$ is the virtual structure sheaf,
${\mathcal L}_i$ is the $i^{th}$ cotangent line bundle, and $\chi$ is
the $K$-theoretic 
push-forward to $\operatorname{Spec}(\cc)$. In algebraic
$K$-theory,
\[
  \chi(M, [F])= \sum_i (-1)^i R^i \pi_* [F]
\]
where $M$ is a variety and $\pi: M \to \operatorname{Spec}(\cc)$ is
the canonical map.

A subring $R\subset H^*(X)$ is {\em self-dual} if the restriction
of the cohomological Poincar\'e pairing to $R$ is nondegenerate.
The
$K$-theoretic Poincar\'e pairing on a nonsingular variety $X$ is: $$
 \la u, v \ra = \chi(X, u \otimes v).$$ 
A subring $R\subset K^*(X)$ is {\em self-dual} if the
restriction of the $K$-theoretic Poincar\'e pairing to $R$ is
nondegenerate.

\begin{theorem} \label{t:2} A reconstruction result from 1-point invariants 
holds in both quantum cohomology and quantum K-theory:
\begin{enumerate}
\item[(i)] Let $R\subset H^*(X)$ be a self-dual subring generated 
by Chern classes of elements of $\Pic(X)$. Let $R^{\perp}$ be the 
orthogonal complement (with respect to the cohomological Poincar\'e pairing). Suppose
\begin{equation*} \label{hcondition}
  (\tau_{k_1}(\gamma_1),\cdots,\tau_{k_{n-1}}(\gamma_{n-1}),
 \tau_{k_{n}}(\xi))_{0,n, \beta} = 0
\end{equation*}
for all $n$-point descendent invariants satisfying
 $\gamma_i \in R$ and $\xi \in R^{\perp}$. 
Then, all $n$-point descendent invariants of classes
of $R$ can be reconstructed from 1-point descendent 
invariants of $R$.
\item[(ii)] Let $R \subset K^*(X)$ be a self-dual subring generated 
by elements of $\Pic(X)$. Let $R^{\perp}$ be the 
orthogonal complement
(with respect to the  Poincar\'e pairing in $K$-theory).
%\footnote{The
%$K$-theoretic Poincar\'e pairing on a nonsingular variety $V$ is: $$
% \la u, v \ra = \chi(V, u \otimes v).$$
%}). 
Suppose 
\[
  (\tau_{k_1}(\gamma_1),\cdots,\tau_{k_n}(\gamma_n),
  \tau_{k_{n+1}}(\xi))^K_{0,n,\beta} =0
\]
for all $n$-point descendent invariants satisfying
$\gamma_i \in R$ and $\xi \in R^{\perp}$. 
Then, all $n$-point descendent invariants of classes
from $R$ can be reconstructed from 1-point descendent 
invariants of $R$. 
\end{enumerate}
\end{theorem}

Part (i) of Theorem 2 is a direct consequence of
Corollary \ref{c:1}, the string equation, and the splitting
axiom of Gromov--Witten theory. 
The self-dual and vanishing conditions on $R$ are required to
control 
the K\"unneth components of the
diagonal arising in the splitting axiom. 
Part (ii) is proven by a parallel argument
in quantum
$K$-theory \cite{yl2}.
We note the subring $R$ need {\em not} be generated
by the entire Picard group for either part of Theorem 2.

In the Gromov--Witten case (i), a similar reconstruction result was proven 
independently by A.~Bertram and H.~Kley in \cite{bk} 
using a very different technique: recursive
relations are found via a residue analysis of the virtual localization
formula of \cite{gp} applied to the graph space of $X$. Our recursive
equations differ from \cite{bk}.
We point out the self-dual condition on $R$ was omitted in the 
Bertram--Kley result \cite{bk} in error.

\subsection{Applications}
\subsubsection{} The 1-point
descendents are the most accessible integrals in quantum
cohomology. Their generating
function, the $J$-function, has been explicitly computed for many
important target varieties (for example, toric varieties and homogeneous 
spaces). The $n$-point descendent invariants, however, remain largely unknown.
Theorem~\ref{t:2} yields a reconstruction of all gravitational
Gromov--Witten invariants from the $J$-function in the case $H^*(X)$ is
generated by $\Pic(X)$. 

For example,
the 1-point invariants of all flag spaces $X$ (associated to simple Lie 
algebras) have been computed by B.~Kim \cite{bkim}. 
As a result, a presentation of the quantum cohomology ring
$QH^*(X)$ can be found. However, 
the 3-point invariants, or structure constants of 
$QH^*(X)$, remain mostly unknown. The 3-point invariants
for the flag space of $A_n$-type have been determined by
Fomin--Gelfand--Postinikov \cite{fgp}.
It is hoped that Theorem \ref{t:1} may help to find a solution for other
flag spaces.
Computations in this direction have been done by H.~Chang and the first author
(in agreement with results of Fomin--Gelfand--Postinikov).
The principal difficultly is to understand the combinatorics
associated to Theorem 2.

\subsubsection{}
The \emph{Quantum Lefschetz Hyperplane Theorem} \cite{yl1}
determines the 1-point descendents of the restricted
classes $i_Y^* H^*(X)$ of a nonsingular very ample divisor
$$i_Y: Y \hookrightarrow X$$
from the 1-point descendents of $X$. 
The following Lemma shows  Theorem~\ref{t:2} may be applied
to the subring $i_Y^*H^*(X) \subset H^*(Y)$.

\begin{lemma}
Let $Y$ be a nonsingular very ample divisor in $X$ determined by the zero locus
of a line bundle $E$.
Assume:
\begin{enumerate}
\item[(i)]
$H^*(X)$ is generated by $Pic(X)$,
\item[(ii)]
$i_{Y*}: H_2(Y) \stackrel{\sim}{\rightarrow} H_2(X)$.
\end{enumerate}
 Consider the ring $R=i_Y^* H^*(X)\subset H^*(Y)$. 
Then, $R$ is self-dual, the vanishing condition for quantum cohomology
in part (i) of Theorem~\ref{t:2} is satisfied, 
and the reconstruction result  holds.
\end{lemma}

\begin{proof}
We first prove $R\subset H^*(Y)$ is self-dual for the
Poincar\'e pairing on $Y$. Equivalently, we will prove, for
all non-zero $\epsilon_Y \in R$, there exists a element
$\delta_Y \in R$ such that
$$\int_Y \epsilon_Y \cup \delta_Y \neq 0.$$
%Let $E \to X$ be the line bundle whose regular zero locus is $Y$
%and 
Let $\epsilon \in H^*(X)$ pull-back to $\epsilon_Y$:
$i_Y^*(\epsilon)=\epsilon_Y.$ 
If 
$$\epsilon \cup c_1(E) \neq 0 \in H^*(X),$$
then there exists $\delta \in H^*(X)$ satisfying:
$$\int_X \epsilon \cup \delta \cup c_1(E) \neq 0$$
as the Poincar\'e pairing on $X$ is nondegenerate.
Let $\delta_Y =i_Y^*(\delta)$. Then,
$$\int_Y \epsilon_Y \cup \delta_Y = \int_X \epsilon \cup \delta \cup c_1(E) \neq 0.$$

If $\epsilon \cup c_1(E) =0\in H^*(X)$, then we will apply the Hard
Lefschetz Theorem (HLT) to prove $\epsilon_Y=0 \in H^*(Y)$.
Let $n$ be the complex dimension of $X$. By HLT applied to $(X, c_1(E))$, we may assume
$$\epsilon \in H^{n-1+k}(X)$$ for $k>0$. Thus,
$\epsilon_Y \in H^{n-1+k}(Y)$. By HLT applied to $(Y, c_1(E_Y))$, there
exists an element $\epsilon'_Y \in H^{n-1-k}(Y)$ satisfying:
$$ \epsilon'_Y \cup c_1(E_Y)^k = \epsilon_Y \in H^*(Y).$$
By the Lefschetz Hyperplane Theorem applied to $Y\subset X$,
$$i_Y^*: H^{n-1-k}(X) \stackrel{\sim}{\rightarrow} H^{n-1-k}(Y).$$
Let $\epsilon'\in H^{n-1-k}(X)$ satisfy $i_Y^*(\epsilon')=\epsilon'_Y$.
We find,
$$ i_Y^*(\epsilon' \cup c_1(E)^k) = \epsilon_Y.$$
As $i_Y^*(\epsilon)=\epsilon_Y$ and $\epsilon \cup c_1(E) =0 \in H^*(X)$, we
find:
$$\epsilon' \cup c_1(E)^{k+1} = 0 \in H^*(X).$$
By HLT applied to $(X, c_1(E))$, we conclude $\epsilon'=0$. The vanishing
$\epsilon_Y=0$ then follows.

Next, we prove the vanishing of Gromov--Witten invariants
required for Theorem 2.
Let $\gamma_i \in H^*(X)$ and $\xi \in R^{\perp}$.
$$
  (\tau_{k_1}(i_Y^* \gamma_1),\cdots,\tau_{k_{n-1}}(i_Y^* \gamma_{n-1}),
 \tau_{k_{n}}(\xi))^Y_{0,n, \beta} 
$$
$$
 =  \int_{[\mbar_{0,n}(Y,\beta)]^{\vir}} \prod_{i=1}^{n-1} \psi_i^{k_i} 
  \ev_i^*(i_Y^*\gamma_i) \ \psi_n^{k_n}  \ev_n^*(\xi).  $$
Since  ${i_Y}_* [\mbar_{0,n}(Y,\beta)]^{\vir}= 
  c_{top}(E_{\beta}) \cap [\mbar_{0,n}(X,\beta)]^{\vir}$ (see \cite{ckl})
 and $$\prod_{i=1}^{n-1} \psi_i^{k_i} \ev_i^*(\gamma_i)\ \psi_n^{k_n}$$
  is a cohomology class pulled-back from $\mbar_{0,n}(X,\beta)$, the
above integral may be rewritten as:
$$
 \int_{Y} i_Y^* ({\ev_n}_* \prod_{i=1}^{n-1} \psi_i^{k_i} \ev_i^*(\gamma_i)\ \psi_n^{k_n} 
  \cap [\mbar_{0,n}(X,\beta)]^{\vir}) \cup  \xi
 =0.$$
\end{proof}

The Quantum Lefschetz Hyperplane Theorem and Theorem~\ref{t:2}
allow the determination of $n$-point descendents of $i_Y^*H^*(X)$ from 
the 1-point descendents of $X$ in this case (see \cite{yl1}, Corollary 1).

The Quantum Lefschetz Hyperplane Theorem and
Lemma 1 also hold when $Y \subset X$ is the nonsingular complete intersection
of very ample divisors. The proof of Lemma 1 for complete intersections is
the same.

\subsubsection{} \label{s:1.3.3} 
Quantum $K$-theory is more difficult than Gromov--Witten theory.
For example, there are no dimension restrictions for the quantum $K$-theoretic
invariants. The genus 0 reconstruction results of Kontsevich-Manin
via the WDVV-equations are less effective in quantum $K$-theory:
the 3-point invariants needed for reconstruction are non-trivial
even in the case of projective space.
Theorem \ref{t:2} provides a new
tool for the study of  quantum $K$-theory.

In case $X=\pp^r$, the 1-point quantum $K$-invariants
have been determined in \cite{yl2} (see Section 2). 
Theorem~\ref{t:2} then allows a recursive computation 
of all the genus 0 $K$-theoretic invariants. The set of invariants includes
(and is essentially equivalent to) the holomorphic Euler characteristics
of the {\em Gromov--Witten subvarieties} of $\mbar_{0,n}(\pp^r,\beta)$.
The Gromov--Witten subvarieties are defined by:
\begin{equation}
\label{sdfds}
\ev_1^{-1}(P_1) \cap \ev_2^{-1}(P_2) \cap \cdots \cap \ev_n^{-1}(P_n) 
\subset \mbar_{0,n}(\pp^r,\beta),
\end{equation}
where $P_1, P_2, \ldots, P_n\subset \pp^r$ are general linear 
subspaces.
By Bertini's Theorem, the intersection (\ref{sdfds}) is a nonsingular
substack. The holomorphic Euler characteristics of
the Gromov--Witten subvarieties specialize to enumerative invariants
for $\pp^r$ when the intersection (\ref{sdfds}) is 0 dimensional.

The $K$-theoretic application was our primary motivation for the
study of the linear relations on the moduli space of maps appearing in 
Theorems 1 and 2.

\subsubsection{} A.~Givental has informed us that our equation \eqref{e:1} 
has a natural interpretation in \emph{symplectic field theory}.

\section{Proofs}

\subsection{Proof of  Theorem~\ref{t:1}}
Consider the moduli space $\mbar_{0,n}(\pp^r, \beta)$.
The curve class is a multiple of the class of a line: $\beta = d[line]$.
If $d=0$, then $n\geq 3$ by the definition of stability.
Equation (\ref{e:1}) is trivial in the $d=0$ case. 
Equation (\ref{e:2}) is easily verified for $\mbar_{0,3}(\pp^r, 0)$. For $d=0$
and $n>3$, the
second equation is obtained by pull-back from the 3-pointed case.
We may therefore assume $d>0$.

It is sufficient to prove 
equations (\ref{e:1}-\ref{e:2}) 
on the 2-pointed moduli space $\mbar_{0,2}(\pp^r, \beta)$ with
marking set $\{i,j\}$.
The equations on the $n$-pointed space $\overline{M}_{0,n}(\pp^r,\beta)$
are then obtained by
pull-back.

Let $B \hookrightarrow \overline{M}_{0,2}(\pp^r,\beta)$ be
a nonsingular curve intersecting the boundary divisors transversely
at their interior points. By the main results of \cite{rp},  equations
(\ref{e:1}-\ref{e:2}) may be established in the Picard
group of $\overline{M}_{0,2}(\pp^r,\beta)$ by proving
the equalities hold after intersecting with all such curves $B$ (actually much
less is needed).

Consider the following fiber square:
\[
 \begin{CD}
  S @>>> C @>f>> \pp^r \\
  @VV{\pi}V @VV{\pi_C}V \\
  B @>>> \mbar_{0,2}(\pp^r,\beta)
 \end{CD}
\]
where $C$ is the universal curve
and $S$ is a nonsingular surface. The morphism 
$\pi$ has sections $s_i$ and $s_j$ induced
from the marked points of $\pi_C$. We find:
\[
 \langle B, \ev_i^* L \rangle =
  \langle s_i, f^* L\rangle, 
\]
\[
  \langle B,\ev_j^* L \rangle =
  \langle  s_j, f^* L\rangle. 
\]
$$ \langle B,\psi_i \rangle = -\langle s_i, s_i \rangle,$$
$$ \langle B, \psi_j \rangle = -\langle s_j, s_j \rangle,$$
where the right sides are all intersection products in $S$.

$S$ is a $\pp^1$-bundle $P$ over $B$ blown-up over points where
$B$ meets the boundary divisors. More precisely, each
reducible fiber of $S$ is a union of two ($-1$)-curves. After a
blow-down of one ($-1$)-curve in each reducible fiber, a $\pp^1$-bundle
$P$ is obtained. Let $P=\pp(V)$
where $V\to B$ is a rank two bundle.
Therefore, 
\[
  0 \to \Pic(P) \to \Pic(S) \to \bigoplus_{b\in Sing} \zz E_b \to 0,
\]
where $Sing \subset B$ in the set of points $b\in B$ where $S_b$ is singular. $E_b$
is the corresponding exceptional divisor of $S$. 
A line bundle $H$ on $S$ is uniquely determined by three sets of
invariants $(J,d, \{e_b\})$, where $H$ is (the pull-back of) an element of
$Pic(B)$, $d$ is the fiber degree of the $\pp^1$-bundle $P$ over $B$, and
$\{e_b\}$ is the set of degrees on the exceptional divisors $E_b$:
$$H = \pi^*(J) \otimes \so_P(d) (- \sum_{b\in Sing} e_b E_b).$$

We may assume $L=\so_{\pp^r}(1)$. Then, $f^* L$ is a line bundle on $S$ of
type
$(J, d, \{d_b \})$ where
$d_b$ is the degree of the map $f$ on 
the exceptional divisors. Similarly,
the sections
$s_i$ and $s_j$ are divisors on $S$  of type
$(J_i, 1, \{\delta_b^i\})$  and
$(J_j, 1, \{\delta_b^j\})$ respectively. Here, $\delta_b^i=1$ or $0$ if
$s_i$ does or does not intersect $E_b$ (and similarly for $\delta_b^j$).

By intersection calculations in $S$, we find:
\begin{equation*} 
 \langle s_i, f^*L \rangle = deg(J)+ d \cdot deg(J_i) + d\cdot
c_1(V)  - \sum_{b\in Sing} d_b \delta_b^i,
\end{equation*}

\begin{equation*}
 \langle s_j, f^*L  \rangle = deg(J)+ d \cdot deg(J_j) + d\cdot
c_1(V)  - \sum_{b\in Sing} d_b \delta_b^j,
\end{equation*}

\begin{equation*}
 -\langle s_i, s_i \rangle = - 2 deg(J_i)- c_1(V)  + \sum_{b \in Sing} \delta_b^i,
\end{equation*}

\begin{equation*}
 -\langle s_j, s_j \rangle = - 2 deg(J_j)- c_1(V)  + \sum_{b \in Sing} \delta_b^j,
\end{equation*}
As $(\pi:S\to B, f: S \to \pp^r, s_i,s_j)$ is
a family of stable maps, the relation
\begin{equation} \label{e:4}
 \langle s_i , s_j\rangle= deg(J_i) + deg(J_j)+ c_1(V) - \sum_{b\in Sing} \delta_b^i \delta_b^j =0 
\end{equation}
is obtained.

Let $Sing(i) \subset Sing$ denote the points $b$ such that
$s_i$ intersects $E_b$ and $s_j$ does not.
Similarly, let $Sing(j) \subset Sing$ denote the 
subset where $s_j$ intersects $E_b$ and $s_i$ does not.

The intersection of equations (\ref{e:1}) and (\ref{e:2})  with
$B$ are easily proven by the above intersection calculations:
$$ \langle f^*L ,s_i \rangle =
 \langle f^*L, s_j \rangle  -d \langle s_j, s_j \rangle 
 - \sum_{b\in Sing(i)} d_b - \sum_{b\in Sing(j)} (d-d_b),$$
$$-\langle s_i, s_i \rangle - \langle s_j, s_j \rangle =
\sum_{b\in Sing(i)} 1 + \sum_{b\in Sing(j)} 1.$$

\noindent
The proof of Theorem 1 is complete.

\subsection{Proof of the Theorem~\ref{t:2}} The result is obtained by an easy 
induction on the number of marked points $n$ and the degree
$\beta$. For simplicity, we assume $Pic(X)={\mathbb Z} H$, $H^*(X)$ is
generated by elements of $Pic(X)$, and $R=H^*(X)$. 
The general argument is identical.

An $n$-point invariant with classes in $H^*(X)$ may be written as
\begin{equation}
\label{dfg}
  \langle \psi^{l_1} H^{k_1}, \cdots, \psi^{l_n} H^{k_n}
   \rangle_{0,n,\beta}.
\end{equation}
Suppose that all $(n-1)$-point invariants and $n$-point
invariants with degree strictly less than $\beta$ are known.

An application of the
first equation of Corollary 1
in case $i=n, j=1$ together with the splitting axiom of Gromov--Witten
theory relates the invariant (\ref{dfg}) to
the invariant 
\begin{equation*}
  \langle \psi^{l_1} H^{k_1+1}, \cdots, \psi^{l_n} H^{k_n-1}
   \rangle_{0,n,\beta}.
\end{equation*}
modulo products of invariants with classes in $H^*(X)$ with either
fewer points or lesser degree. The self-dual and vanishing conditions on $R$ in part (i)
of Theorem 2 are required to kill the diagonal splittings 
not consisting of classes of $R$ --- of course these conditions are
trivial in case $R=H^*(X)$.
After repeating the procedure, we may assume $k_n=0$.

Similarly, applications of the second equation of Corollary 1
allow a reduction of $l_n$ to 0 (modulo known invariants).

Once $l_n=0$ and $k_n=0$, then the $n$-point invariant may be reduced 
to $(n-1)$-point
invariants by the string equation. This completes the induction step.

The proof of Theorem 1 in quantum cohomology is complete. The
argument for quantum $K$-theory is identical. 
One simply replaces the divisor classes (and their products) by the 
$K$-products of the corresponding line bundles.
The splitting axiom of quantum $K$-theory is slightly more complicated
(see \cite{yl2}).

\section{Examples}

\subsection{Gromov--Witten invariants of $\pp^2$}
 Kontsevich's formula for the genus 0
Gromov--Witten invariants of $\pp^2$ is derived here from Theorem 2.

The cohomology ring $H^*(\pp^2)= \qq[H]/(H^3)$ has a linear basis 
$$1, H, H^2.$$ 
By the fundamental class and divisor axioms,
only invariants of the form
$$N_d=(H^2,\cdots,H^2)_{0,3d-1,d}$$ need be computed.
As there is a unique line
through two distinct points in $\pp^2$, we see $N_1=1$.

We may reformulate equation (1) in the
following form, which is better suited for computations without descendents.

\begin{proposition}
Let $n \ge 3$. Let $i,j,k$  be distinct markings. Then
\begin{equation} \label{e:3}
 \ev_i^*(H)= \ev_j^*(H) + \sum_{d_1+d_2=d} 
 \left( d_2   D_{ik,d_1| j,d_2} -d_1 D_{i,d_1|jk,d_2} \right).
\end{equation}  
\end{proposition}  

\begin{proof}
When $n\ge 3$, $\psi_j=D_{j,ik}$. Then, equation (1) easily implies 
the Proposition.
\end{proof}

Consider the following cohomology class in 
$H^*(\overline{M}_{0,{3d-1}}(\pp^2,d))$:
\begin{equation}
\label{gegg}
\ev_1^*(H^2) \cdots \ev_{3d-2}^*(H^2) \ev_{3d-1}^*(H).
\end{equation}
Let $i=3d-1$ and $j=1,k=2$. 
After intersecting \eqref{e:3} with (\ref{gegg}) and applying
the splitting axioms, we obtain: 
\[
  {N_d} =
  \sum_{d_1+d_2=d, d_i>0} N_{d_1} N_{d_2} \left( 
     d_1^2 d_2^2 \binom{3d-4}{3d_1-2}
    - {d^3_1 d_2} \binom{3d-4}{3d_1-1} \right).
\]

\subsection{Quantum $K$-invariants of $\pp^1$}
We explain here the computation of quantum $K$-invariants of
$\pp^1$ in genus 0 using Theorem~\ref{t:1}.

Define the $K$-theoretic $J$-function of $\pp^r$ to be
\[
 J^K_{\pp^r}(Q,q):= \sum_{d=0}^{\infty} Q^d {\ev_d}_* (\frac{1}{1-q \mathcal{L}})
\]
where ${\ev_d}_*: K(\mbar_{0,1}(\pp^r,d)) \to K(\pp^r)$ is the $K$-theoretic
push-forward. The virtual structure sheaf in this case is just
the ordinary structure sheaf as $\mbar_{0,n}(\pp^r,d)$ is a smooth stack.
The following result is proven in \cite{yl2}.

\begin{proposition} \label{p:2}
\[
 J^K_{\pp^r}(Q,q) = \sum_d \frac{Q^d}{\prod_{m=1}^d (1-q^m H)^{r+1}},
\]
where $H=\so(1)$ is the hyperplane bundle in $\pp^r$.
\end{proposition}

We will specialize now to $\pp^1$. By the classical result $$K^*(\pp^1)= 
\frac{\qq[H]}{(H-1)^2},$$  $K^*(\pp^1)$ is a two dimensional 
$\qq$-vector space with basis $$e_0=\so, \ e_1=H-\so.$$ 
The $K$-theoretic Poincar\'e metric is 
\[
  (g_{ij})= \left( 
  \begin{matrix}
  1 & 1 \\
  1 & 0
  \end{matrix}
  \right)
\]
with inverse matrix
\[
(g^{ij})= \left( 
  \begin{matrix}
  0 & 1 \\
  1 & -1
  \end{matrix}
  \right)
\]

From Proposition~\ref{p:2}, we obtain the 1-point quantum 
$K$-invariants,
 $$(\dfrac{\gamma}{1-q\cl})^K_{g,n,d}=
\sum_k q^k (\tau_k(\gamma))^K_{g,n,d},$$

\begin{eqnarray*}
  (\frac{e_1}{1-q \cl})^K_{0,1,1} & = & 
    1+ 2 q + 3 q^2 + \cdots, \\ 
 (\frac{e_1}{1-q \cl})^K_{0,1,2}
 &  = &  1 + 2q + 5 q^2 + \cdots .
\end{eqnarray*}

Let $\gamma_i\in K^*(\pp^1)$. The {\em $K$-theoretic fundamental class equation} is:
\begin{equation} \label{e:8}
  (\gamma_1\cdots,\gamma_{n-1},e_0)^K_{0,n,d}
  =(\gamma_1\cdots,\gamma_{n-1})^K_{0,n-1,d},
\end{equation}
as the fiber of $\mbar_{0,n}(\pp^1,d) \to \mbar_{0,n-1}(\pp^1,d)$ is
rational.

The string equation is obtained from the geometry of the
morphism 
$\pi: \mbar_{0,n}(X,\beta) \to \mbar_{0,n-1}(X,\beta)$ 
forgetting the last point.
The \emph{$K$-theoretic string equation} \cite{yl2} takes the following form: 
\begin{equation} \label{e:string}
 R^0 \pi_*(\otimes_{i=1}^{n-1} \cl_i^{\otimes k_i}) 
  = \otimes_{i=1}^{n-1} \cl^{\otimes k_i} \ \otimes \ (\so + \sum_{i=1}^{n-1}\sum_{k=1}^{k_i} \cl_i^{-k}),
\end{equation}
$$
 R^1 \pi_* (\otimes_{i=1}^{n-1} \cl^{\otimes k_i})= 0.
$$

We will illustrate the computational scheme of $n$-point quantum 
$K$-invariants for $\pp^1$ by calculating  
 $(e_1,e_1)^K_{0,2,2}$. 
We start with by rewriting the invariant:
 $$ ({e_1}, e_1)^K_{0,2,2} 
  =({e_1}, H)^K_{0,2,2}-(e_1, e_0)^K_{0,2,2}.$$
By equations (1) and \eqref{e:8} we find:
 $$({e_1}, e_1)^K_{0,2,2}= (\cl_1^2 e_1 H,e_0)^K_{0,2,2} - 
   \chi \left(D_{1,d=1|2,d=1}, \cl_1^2\ev_1^* (e_1 H)\right)-(e_1)^K_{0,1,2}. $$
We may use the relation $e_1 H=e_1$ and the $1$-point
evaluations. Together with the string equation 
and splitting axiom, we find:
\begin{equation}
\label{fdg}
 (e_1,e_1)^K_{0,2,2} =
   7-
   ( \cl_1^2e_1, e_a)^K_{0,2,1} g^{ab} (e_b, e_0)^K_{0,2,1}.
\end{equation}
% $$= 7 
%   - ({ e_1}, e_0)^K_{0,2,1} (e_1)^K_{0,1,1} 
%    - ({ e_1}, e_1)^K_{0,2,1} (e_0)^K_{0,1,1} 
%   -({  e_1}, e_1)^K_{0,2,1} (-1) (e_1)^K_{0,1,1}$$
The following invariants are easy to compute by equation (1) and the string equation:
\[
  (e_1,e_1)^K_{0,2,1}=1, \qquad (\cl_1^2e_1, e_1)^K_{0,2,1}=4.
\]
Substitution in equation (\ref{fdg}) yields:
$$(e_1,e_1)^K_{0,2,2}=1.$$
A similar, but much longer, computation shows
 $(e_1,e_1,e_1)^K_{0,3,2}=1$.

We conclude with two remarks about the quantum $K$-theory of
$\pp^1$ and the rationality of the Gromov--Witten subvarieties of
$\overline{M}_{0,n}(\pp^1,d)$.

\begin{enumerate}
\item[(i)] By 
Proposition~\ref{p:2} and the fundamental class equation, we see:
\[ (e_0,e_0,\cdots,e_0)^K_{0,n,d}=1. \]
This result may also be deduced from the rationality of
the moduli space $\mbar_{0,n}(\pp^r,d)$ proven in
 \cite{kp}. 
\item[(ii)]
By the  exact sequence,
\begin{equation*}
  0 \to \so \to \so_{\pp^1}(1) \to \so_p \to 0,
\end{equation*}
where $p$ is a point in $\pp^1$, we find  $e_1=[\so_p]$. Hence,
\begin{equation*}
  \left(e_1, \cdots, e_1 \right)^K_{0,n,d}
  =\chi \left(\cap_{i=1}^n \ev_{i}^{-1} (p_i)\right)
\end{equation*}
where $\ev_i^{-1}(p_i)$ and their intersections are Gromov--Witten subvarieties 
(see \S~\ref{s:1.3.3}). 
For small pairs $(n,d)$,  the space $\cap_{i=1}^n \ev_{i,d}^{-1} (p_i)$
is also rational. For example, rationality certainly holds for the
cases $(2,2)$ and $(3,2)$ discussed in the above computations.
 
It is interesting to ask which Gromov--Witten subvarieties of
$\overline{M}_{0,n}(\pp^1,d)$ are rational. For $\pp^2$, irrational 
Gromov--Witten subvarieties have been found in \cite{pan}.
\end{enumerate}


\begin{thebibliography}{30}

\bibitem{bk}
A.~Bertram and H.~Kley, \emph{New recursions for genus-zero Gromov--Witten 
invariants}, math.AG/0007082.

\bibitem{ckl}
D.~A.~Cox, S.~Katz and Y.-P.~Lee,
\emph{Virtual Fundamental Classes of Zero Loci}, to appear in
\emph{Enumerative geometry in physics}, Contemporary Mathematics, AMS. 
math.AG/0006116.

\bibitem{fp}
W.~Fulton and R.~Pandharipande, \emph{Notes on stable maps and quantum 
cohomology}, Algebraic geometry---Santa Cruz 1995, Proc. Symp. Pure. Math.
62. Part 2, (1997) 45-96.

\bibitem{fgp}
S.~Fomin, S.~Gelfand, and A.~Postnikov, \emph{Quantum Schubert polynomials},
J. Amer. Math. Soc. \textbf{10} (1997), no. 3, 565--596. 

\bibitem{gp}
T.~Graber and R.~Pandharipande, \emph{Localization of virtual classes}, 
Invent. Math. \textbf{135} (1999), no. 2, 487--518.

\bibitem{bkim}
B.~Kim, \emph{Quantum cohomology of flag manifolds $G/B$ and quantum Toda 
lattices}, Ann. of Math. (2) \textbf{149} (1999), no. 1, 129--148.

\bibitem{kp}
B.Kim and R.~Pandharipande, \emph{The connectedness of the moduli space of maps 
to homogeneous spaces},  math.AG/0003168.

\bibitem{yl1}
Y.-P.~Lee, \emph{Quantum Lefschetz hyperplane theorem}, to appear in Invent.
Math.

\bibitem{yl2} 
Y.-P.~Lee, \emph{Quantum $K$-theory I: foundation}, preprint. 
\emph{Quantum $K$-theory II}, in preparation.

\bibitem{rp} 
R.~Pandharipande, \emph{Intersections of $\mathbf{Q}$-divisors on 
Kontsevich's moduli space $\overline{M}_{0,n}(\mathbf{P}^r,d)$ 
and enumerative geometry}. Trans. Amer. Math. Soc. \textbf{351} (1999), 
no. 4, 1481--1505 

\bibitem{pan} R.~Pandharipande,
{\em The canonical class of $\overline{M}\sb {0,n}(\mathbf{P}^r,d)$ and 
enumerative geometry},
Internat. Math. Res. Notices (1997), no. 4, 173--186. 
\end{thebibliography}
\end{document}